# Is there a duality in the classical acceptance of non-constructive, foundational, concepts as axiomatic?

Bhupinder Singh Anand

We consider a philosophical question that is implicit in Bringsjord's paper [Br93]: Is there a duality in the classical acceptance of non-constructive, foundational, concepts as axiomatic?

**1. Mendelson's thesis**

We note that Mendelson [Me90] is quoted there as saying (*parenthetical qualifications added*):

(*i*)   Here is the main conclusion I wish to draw: it is completely unwarranted to say that CT is unprovable just because it states an equivalence between a vague, imprecise notion (effectively computable function) and a precise mathematical notion (partial-recursive function).

(*ii*)   The concepts and assumptions that support the notion of partial-recursive function are, in an essential way, no less vague and imprecise (*non-constructive, and intuitionistically objectionable*) than the notion of effectively computable function; the former are just more familiar and are part of a respectable theory with connections to other parts of logic and mathematics. (The notion of effectively computable function could have been incorporated into an axiomatic presentation of classical mathematics, but the acceptance of CT made this unnecessary.) ... Functions are defined in terms of sets, but the concept of set is no clearer (*not more non-constructive, and intuitionistically objectionable*) than that of function and a foundation of mathematics can be based on a theory using function as primitive notion instead of set. Tarski's definition of truth is formulated in set-theoretic terms, but the notion of set is no clearer (*not more non-constructive, and intuitionistically objectionable*), than that of truth. The model-theoretic definition of logical validity is based ultimately on set theory, the foundations of which are no clearer (*not more non-constructive, and intuitionistically objectionable*) than our intuitive (*non-constructive, and intuitionistically objectionable*) understanding of logical validity.

(*iii*) The notion of Turing-computable function is no clearer (*not more non-constructive, and intuitionistically objectionable*) than, nor more mathematically useful (foundationally speaking) than, the notion of an effectively computable function,

where:

(*a*) The Church-Turing Thesis, CT, is formulated as: A function is effectively computable if and only if it's Turing-computable.

(*b*) An effectively computable function is defined to be the computing of a function by an algorithm.

(*c*) The classical notion of an algorithm is expressed by Mendelson as: an effective and completely specified procedure for solving a whole class of problems. ... An algorithm does not require ingenuity; its application is prescribed in advance and does not depend upon any empirical or random factors.

and, where Bringsjord paraphrases (*iii*) as:

(*iv*)  The notion of a formally defined program for guiding the operation of a TM is no clearer than, nor more mathematically useful (foundationally speaking) than, the notion of an algorithm.

adding that:

(*v*)  This proposition, it would then seem, is the very heart of the matter. If (*iv*) is true then Mendelson has made his case; if this proposition is false, then his case is doomed, since we can chain back by modus tollens and negate (*iii*).

## 2. The concept of "constructive, and intuitionistically unobjectionable"

Now, prima facie, any formalisation of a "vague and imprecise", "intuitive" concept, say C, would normally be intended to capture the concept C, both faithfully and completely, within a constructive, and intuitionistically unobjectionable[1], language L.

Clearly, we could disprove the thesis - that C and its formalisation are interchangeable, hence equivalent - by showing that there is a constructive aspect of C that is formalisable in a constructive language L', but that such formalisation cannot be assumed expressible in L without introducing inconsistency.

However, equally clearly, there can be no way of proving the equivalence, as this would contradict the premise that the concept is "vague and imprecise", hence essentially open-ended in a non-definable way, and so non-formalisable.

---

[1] The terms "constructive" and "constructive, and intuitionistically unobjectionable" are used synonymously both in their familiar linguistic sense, and in a mathematically precise sense. Mathematically, we term a concept as "constructive, and intuitionistically unobjectionable" if, and only if, it can be defined in terms of pre-existing concepts without inviting inconsistency. Otherwise, we understand it to mean unambiguously verifiable, by some "effective method", within some finite, well-defined, language or meta-language. It may also be taken to correspond, broadly, to the concept of "constructive, and intuitionistically unobjectionable" in the sense apparently intended by Gödel in his seminal 1931 paper [Go31a].

Obviously, Mendelson's assertion that there is no justification for claiming Church's Thesis as unprovable must, therefore, rely on an interpretation that differs significantly from the above; for instance, his concept of provability may appeal to the axiomatic acceptability of "vague and imprecise" concepts - as suggested by his remarks.

Now, we note that all the examples cited by Mendelson involve the decidability (computability) of an infinitude of meta-mathematical instances, where the distinction between the constructive meta-assertion - that any given instance is individually decidable (computable) - and the non-constructive meta-assertion - that all the instances are jointly decidable (computable) uniformly - is not addressed explicitly. However, (*a*), (*b*) and (*c*) appear to suggest that Mendelson's remarks relate implicitly to non-constructive meta-assertions.

Perhaps the real issue, then, is the one that emerges if we replace Mendelson's use of implicitly open-ended concepts such as "vague and imprecise", and "intuitive", by the more meta-mathematically meaningful concept of "non-constructive, and intuitionistically objectionable", as indicated parenthetically.

The essence of Mendelson meta-assertion (*iii*), then, appears to be that, if the classically accepted definitions of foundational concepts such as "partial recursive function", "function", "Tarskian truth" etc. are also non-constructive, and intuitionistically objectionable, then replacing one non-constructive concept by another may be psychologically unappealing, but it should be meta-mathematically valid and acceptable.

### 3. The duality

Clearly, meta-assertion (*iii*) would stand refuted by a non-algorithmic effective method that is constructive. However, if it is explicitly - and, as suggested by the nature of the arguments in Bringsjord's paper, widely - accepted at the outset that any effective method is necessarily algorithmic (i.e. uniform as stated in (*d*) below), then any counter-argument to CT can, prima facie, only offer non-algorithmic methods that may, paradoxically, be effective intuitively in a non-constructive, and intuitionistically objectionable, way only!

Recognition of this dilemma is implicit in the admission that the various arguments, as presented by Bringsjord in the case against Church's Thesis - including his narrational case - are open to reasonable, but inconclusive, refutations. Nevertheless, if we accept Mendelson's thesis that the inter-changeability of non-constructive concepts is valid in the foundations of mathematics, then Bringsjord's case against Church's Thesis, since it is based similarly on non-constructive concepts, should also be considered conclusive classically (even though it cannot, prima facie, be considered constructively conclusive in an intuitionistically unobjectionable way). There is, thus, an apparent duality in the - seemingly extra-logical - decision as to whether an argument based on non-constructive concepts may be accepted as classically conclusive or not.

That this duality may originate in the very issues raised in Mendelson's remarks - concerning the non-constructive roots of foundational concepts that are classically accepted as mathematically sound - is seen if we note that these issues may be more significant than is, prima facie, apparent.

**4. Definition of a formal mathematical object, and consequences**

Thus, if we define a formal mathematical object as any symbol for an individual letter, function letter or a predicate letter that can be introduced through definition into a formal theory without inviting inconsistency, then it is simply shown[2] that unrestricted, non-constructive, definitions of non-constructive, foundational, set-theoretic concepts - such as "mapping", "function", "recursively enumerable set", etc. - in terms of constructive number-theoretic concepts - such as recursive number-theoretic functions and relations - do not always correspond to formal mathematical objects.

---

[2] See Meta-lemma 2, and its corollaries, in Anand [An02].

In other words, the above assumption - that every such definition corresponds to a formal mathematical object - introduces a formal inconsistency into standard Peano Arithmetic and, ipso facto, into any Axiomatic Set Theory that models standard PA (loosely speaking, this may be viewed as a constructive arithmetical parallel to Russell's non-constructive impredicative set).

Since it can also be argued[3] that the non-constructive element in Tarski's definitions of "satisfiability" and "truth", and in Church's Thesis, originate in a common, but removable, ambiguity in the interpretation of an effective method, perhaps it is worth considering whether Bringsjord's acceptance of the assumption (*d*) - that every constructive effective method is necessarily algorithmic, in the sense of being a *uniform* procedure - is mathematically necessary, or even whether it is at all intuitively tenable.

   (*d*)   … *uniform* procedure, a property usually taken to be a necessary condition for a procedure to qualify as effective.

Thus, we may argue that we can explicitly, and constructively, define a non-algorithmic, effective method as one that, in any given, individual case, is individually effective if, and only if, it terminates finitely, with a conclusive result; and an algorithmic effective method as one that is uniformly effective if, and only if, it terminates finitely, with a conclusive result, in any given, individual case.[4]

---

[3] See §5 of Anand [An02].

[4] We note that the possibility of a distinction between the interpreted number-theoretic meta-assertions, "For any given natural number *x*, *F*(*x*) is true" and "*F*(*x*) is true for all natural numbers *x*", is not evident unless these are expressed symbolically as, "(A*x*)(*F*(*x*) is true)" and "(A*x*)*F*(*x*) is true", respectively. The issue, then, is whether the distinction can be given any mathematical significance. For instance, under a

## 5. Bringsjord's case against CT

(1) Apropos the specific arguments against CT, it would seem, prima facie, that an individually effective - even if not obviously constructive - method could be implicit in the following argument considered by Bringsjord:

> (*e*) Assume for the sake of argument that all human cognition consists in the execution of effective processes (in brains, perhaps). It would then follow by CT that such processes are Turing-computable, i.e., that computationalism is true. However, if computationalism is false, while there remains incontrovertible evidence that human cognition consists in the execution of effective processes, CT is overthrown.

Assuming computationalism is false, the issue in this argument would, then, be whether there is a constructive, and adequate, expression of human cognition in terms of individually effective methods.

(2) An appeal to such an individually effective method may, in fact, be implicit in Bringsjord's consideration (*italicised comments in parentheses added*) of the predicate *H*, defined by:

$H(P, i)$ iff $(En)S(P, i, n)$

where the predicate $S(P, u, n)$ holds if, and only if, TM *M*, running program *P* on input *u*, halts in exactly *n* steps ($= M_P : u =>_n$ halt).

Bringsjord defines *S* as **totally** (*and, implicitly, uniformly*) **computable** in the sense that, given some triple $(P, u, n)$, there is some (*uniform*) program *P*\* which, running

---

constructive formulation of Tarski's definitions, we may qualify the latter by saying that it can be meaningfully asserted as a totality only if "*F*" is a mathematical object.

on some TM $M^*$, can infallibly give us a verdict, **Y** ("yes") or **N** ("no"), for whether or not $S$ is true of this triple.

He then notes that, since the ability to (*uniformly*) determine, for a pair ($P$, $i$), whether or not $H$ is true of it, is equivalent to solving the full halting problem, $H$ is not totally computable. However, he also notes that there is a program[5] (*implicitly non-uniform, and so, possibly, effective individually*) which, when asked whether or not some TM $M$ run by $P$ on $u$ halts, will produce **Y** iff $M_P : u =>_n$ halt. For this reason $H$ is declared **partially** (*implicitly, individually*) **computable**.

(3) More explicitly, Bringsjord remarks that Kalmár's refutation of CT is classically inconclusive mainly because it does not admit any uniform effective method, but appeals to the existence of an infinitude of individually effective methods. However, Kalmár's argument can be strengthened if we note, firstly, that we can extend the definition of classical Turing machines to include meta-routines that constructively self-terminate the computational process of a Turing machine whenever an instantaneous tape description repeats itself.

Secondly, a constructive formulation of Tarski's definitions of "satisfiability" and "truth", and of a Uniform Church Thesis, then implies[6] that every partial recursive number-theoretic function can be constructively, and uniquely, extended as a total function that is not classically Turing-computable - thus disproving the classical Church-Turing Thesis conditionally, and effectively solving the Halting problem, also conditionally.

---

[5] An effective, conditional, solution to the Halting problem is given in Anand ([An02], Corollary 14.1).

[6] See Meta-lemma 14 in Anand [An02].

Excepting, that it always calculates $g(n)$ constructively, even in the absence of a uniform procedure, within a fixed postulate system, the reasoning used in the above proof is, essentially, the same as the argument in Kalmár [Ka59], reproduced below from Bringsjord's paper:

> First, he draws our attention to a function $g$ that isn't Turing-computable, given that $f$ is[7]:
>
> $$g(x) = \mu_y(f(x, y) = 0) = \{\text{the least } y \text{ such that } f(x, y) = 0 \text{ if } y \text{ exists; and } 0 \text{ if there is no such } y\}$$
>
> Kalmár proceeds to point out that for any $n$ in **N** for which a natural number $y$ with $f(n, y) = 0$ exists, "an obvious method for the calculation of the least such $y$ ... can be given," namely, calculate in succession the values $f(n, 0), f(n, 1), f(n, 2), ...$ (which, by hypothesis, is something a computist or TM can do) until we hit a natural number $m$ such that $f(n, m) = 0$, and set $y = m$.
>
>> On the other hand, for any natural number $n$ for which we can prove, not in the frame of some fixed postulate system but by means of arbitrary - of course, correct - arguments that no natural number $y$ with $f(n, y) = 0$ exists, we have also a method to calculate the value $g(n)$ in a finite number of steps.[8]
>
> Kalmár goes on to argue as follows. The definition of $g$ itself implies the *tertium non datur*, and from it and CT we can infer the existence of a natural number $p$ which is such that
>
> (\*)   there is no natural number $y$ such that $f(p, y) = 0$; and (\*\*) this cannot be proved by any correct means.

---

[7] Bringsjord notes that the original proof can be found on page 741 of Kleene [Kl36].

[8] Quoted by Bringsjord from Kalmár [Ka59].

Kalmár claims that (\*) and (\*\*) are very strange, and that therefore CT is at the very least implausible.

(*Updated: Saturday 10${}^{th}$ May 2003 1:36:52 AM by re@alixcomsi.com*)